\documentclass[11pt]{amsart}
\usepackage{amsmath,amssymb, amsthm}
\usepackage{setspace}
\usepackage[all]{xy}
\newcommand\Seq{\xymatrix{
1\ar[r]&D\ar[r]&{\Delta}\ar[r]&{\Gamma}\ar[r]&1\\
}}

\newcommand\fn{\mathfrak n}

\newcommand\inv{^{-1}}
\newcommand\eg{E(\gamma)}
\newcommand\RB{{\mathbb R}}

\DeclareMathOperator{\Lie}{Lie}

\DeclareMathOperator{\Aut}{Aut}
\DeclareMathOperator{\Out}{Out}

\newtheorem{theorem}{Theorem}[section]

\newtheorem{lemma}[theorem]{Lemma}
\newtheorem{corollary}[theorem]{Corollary}
\newtheorem{defn}[theorem]{Definition}

\newtheorem{conjecture}[theorem]{Conjecture}
\newtheorem{question}[theorem]{Question}

\begin{document}

\title{Continuous Quotients for Lattice Actions on Compact Spaces}
\author{David Fisher and Kevin Whyte}

\begin{abstract}
Let $\Gamma<SL_n({\mathbb Z})$ be a subgroup of finite index, where 
$n{\geq}5$. 
Suppose $\Gamma$ acts continuously on a manifold $M$, where 
$\pi_1(M)={\mathbb Z}^n$, 
preserving a measure that is positive on open sets.
Further assume that the induced $\Gamma$ action on $H^1(M)$ is 
non-trivial.  We show there exists a finite index subgroup 
$\Gamma'<\Gamma$ and a $\Gamma'$ equivariant continuous map 
$\psi:M{\rightarrow}{\mathbb T}^n$ that induces an isomorphism on 
fundamental group.
        
We prove more general results providing continuous quotients in cases 
where $\pi_1(M)$ surjects onto a 
finitely generated torsion free nilpotent group.  We also give some 
new examples of manifolds with 
$\Gamma$ actions.
\end{abstract}

\maketitle

\section{{\bf Introduction}}
\label{section:introduction}
Let $G$ be a semisimple Lie group with $\mathbb R$-rank$(G){\geq}2$, 
and $\Gamma<G$ a lattice.  In this paper 
we seek to examine the relationship between the topology and dynamics 
of measure preserving actions of 
$\Gamma$ on a compact manifold $M$.  More precisely, we seek 
relations between the fundamental group of 
$M$ and the structure of both $M$ and the $\Gamma$ action on $M$.  
The simplest version of our main theorem is:

\begin{corollary}
\label{corollary:torus}
Let $\Gamma<SL_n(Z), n{\geq}3$ be a subgroup of finite index.  
Suppose $\Gamma$
acts on a compact manifold $M$ preserving a measure that is positive 
on open sets.
Assume ${\rho}:\pi_1(M){\rightarrow}{\mathbb Z}^n$ is $\Gamma$ 
equivariant, that the action of 
$\Gamma$ on ${\mathbb Z}^n$ is given by the standard representation 
of $SL_n(\RB)$, and that the $\Gamma$
action lifts to $\ker{\rho}$. (The lifting condition is automatic 
provided $n{\geq}5$).  Then
there is a finite index subgroup $\Gamma'<{\Gamma}$ and a $\Gamma'$ 
equivariant map 
$\psi:M{\rightarrow}{\mathbb T}^n$ which induces the map $\rho$ on 
fundamental groups.
\end{corollary}

The actual result applies more generally to certain kinds of actions 
(described precisely below) of 
lattices on manifolds whose fundamental group surjects onto a torsion 
free finitely generated 
nilpotent group.  In all cases we produce a continuous map from 
our action to an algebraically defined action on a nilmanifold.  In 
particular, the theorem applies to
the surgery examples of Katok and Lewis, and the quotient recaptures 
the torus action 
on which the surgeries were performed \cite{KL}.  For some of the 
examples with more complicated fundamental
group it will be necessary to pass to a finite cover before our 
theorem applies.  For all of these examples, the 
map to the torus simply collapses certain invariant submanifolds 
whose fundamental groups map to $0$ in 
${\mathbb Z}^n$.  In section \ref{section:a&e}, we prove a purely 
topological result that shows that this is essentially what all 
such maps must look like.

This work is related to work of the first author, who shows in 
\cite{F} that under different hypothesis
there are measurable maps to the same standard examples.  The 
argument of \cite{F} also shows that 
fundamental groups of manifolds with measure preserving $\Gamma$ 
actions are of arithmetic type under 
fairly mild hypotheses relating the dynamics  to the topology, unless 
the fundamental group admits no linear 
representation.  Further work of the first author with Zimmer 
\cite{FZ} produces linear representations of the 
fundamental group under stronger geometric assumptions.  Combining 
these results gives that either the action
on fundamental group is trivial or that the fundamental group 
contains a direct factor that is finitely 
generated and nilpotent on which $\Gamma$ acts nontrivially.  These 
results give some evidence that our 
assumptions are not atypical, at least in the case that 
$\Gamma{\rightarrow}(\Out(\pi_{1}(M)))$ is nontrivial.

We also discuss several examples that show that for our general 
result one cannot expect 
to improve the regularity of the quotients without stronger 
assumptions.  Note that our theorems apply
not just to manifolds but to any topological spaces satisfying 
standard covering space theory, i.e. any connected,
locally path-connected, semi-locally 1-connected, locally compact, 
separable, metrizable space with 
finitely generated fundamental group. 

\section{{\bf Preliminaries}}
\label{section:preliminaries}

We assemble here some basic facts from topology and dynamics that 
will be used in the proof.

Suppose a finitely generated group $\Gamma$ acts continuously on a 
manifold $M$.  Let $\tilde M$ be any 
cover of $M$ and $D$ and $\Delta$ be the deck group of $\tilde M$ over 
$M$ and the group of lifts of the 
$\Gamma$ action, respectively.  There is an exact sequence:
$$\Seq$$
and the $\Gamma$ action lifting to $\tilde M$ is equivalent to this 
exact sequence splitting.  In 
particular, elementary group cohomology shows that the sequence 
splits if $H^2({\Gamma}, D)=0$ and the map 
$\Gamma{\rightarrow}\Out(D)$ given by the sequence lifts to a map 
$\Gamma{\rightarrow}\Aut(D)$.

Suppose $\Gamma$, a group, acts on a measure space $M$ preserving a 
finite measure.  A cocycle is a 
measurable map $\alpha:{\Gamma}{\times}M{\rightarrow}H$ where $H$ is 
a group and $\alpha$ satisfies 
the equation $\alpha(g_1g_2, m)=\alpha(g_1,g_2m)\alpha(g_2, m)$ for 
all $g_1,g_2{\in}{\Gamma}$ and 
all $m{\in}M$.
  Two cocycles $\alpha$ and $\beta$ are cohomologous if there is a 
measurable map $\phi:M{\rightarrow}H$ 
such that $\alpha(g, m)={\phi(gm)\inv}\beta(g,m)\phi(m)$ for all 
$g{\in}{\Gamma}$ and almost all $m{\in}M$.
  Now assume that $\Gamma<G$ is a lattice where $G$ is a semisimple 
Lie group, that the $\Gamma$ action on
 $M$ is ergodic,  and that $H$ is an algebraic group.  By results of 
Zimmer, for any cocycle 
$\alpha:{\Gamma}{\times}M{\rightarrow}H$, there is a minimal subgroup 
$L<H$, unique up to conjugacy, 
such that $\alpha$ is cohomologous to a cocycle taking values in 
$L$.  The group $L$ is called the 
algebraic hull of the cocycle \cite{Z1}.

We will also need the following definition:

\begin{defn}
A representation $\sigma:{\Gamma}{\rightarrow}{GL_n(\mathbb R)}$ is 
weakly hyperbolic if there is no invariant subspace $V<{\mathbb R}^n$ 
where all the eigenvalues of all elements of $\Gamma$ have modulus 
one.
\end{defn}

\section{{\bf Main Theorem and Proof}}
\label{section:maintheorem}

Suppose $\Gamma$, a finitely generated group,  acts on a compact 
manifold $M$ and there is a surjection 
$\rho:\pi_1(M){\rightarrow}{\Lambda}$, where ${\Lambda}$ is a 
finitely generated torsion free nilpotent 
group.

If the action of $\Gamma$ lifts to the cover of $M$ corresponding to 
$\ker{\rho}$, then we have a map 
$\sigma:{\Gamma}{\rightarrow}\Aut(\Lambda)$.  By a theorem of Malcev, 
this gives a map 
$\sigma:\Gamma{\rightarrow}\Aut(N)$ where $N$ is a simply connected 
nilpotent Lie group, containing 
$\Gamma$ as a lattice.  We can view $\sigma$ as a representation 
$\sigma:{\Gamma}{\rightarrow}{GL(\fn)}$ 
where $\fn=\Lie(N)$.

\begin{defn}
Under the conditions described above, we say that the $\Gamma$ action 
on $M$ is $\pi_1$-hyperbolic if 
$\sigma:{\Gamma}{\rightarrow}GL(\fn)$ is weakly hyperbolic.
\end{defn}

Note that our definition of $\pi_1$ hyperbolic includes the 
assumption that the $\Gamma$ action lifts to 
the appropriate cover of $M$.

\begin{theorem}
\label{theorem:main}
Let ${\Gamma}<G$ be an irreducible lattice, where $G$ is a semisimple 
Lie group with all simple factors of ${\mathbb R}$-rank$(G){\geq}2$.  
Suppose that $\Gamma$ acts on a
compact manifold $M$ preserving a measure that is positive on open 
sets.  Assume there is
a surjection  $\rho:{\pi_1(M)}{\rightarrow}{\Lambda}$ where
$\Lambda$ is a finitely generated torsion free nilpotent group and
that the $\Gamma$ action is ${\pi_1}$-hyperbolic.  Then there is a
finite index subgroup $\Gamma'<{\Gamma}$ and a  $\Gamma'$ equivariant 
map
$\psi:M{\rightarrow}N/{\Lambda}$ where the action on $N/{\Lambda}$ is 
defined by extending the
action of $\Gamma'$ on $\Lambda$ to $N$.  Furthermore, the map 
$\psi_*$ on fundamental group is equal to 
$\rho$ above.
\end{theorem}
 
\noindent
In the case of $M$ having abelian fundamental group, this gives 
corollary \ref{corollary:torus} from
the introduction.  We need only explain the remark that the lifting 
of the action is automatic in the case
where $n{\geq}5$.  To see this, we need only see that 
$H^2(\Gamma, {\mathbb Z}^n)=0$, since $\Aut({\mathbb 
Z}^n)=\Out({\mathbb Z}^n)$.  By theorem $4.4$ of \cite{B},
we know that $H^2(\Gamma, {\mathbb R}^n)=0$ .  By looking at the long 
exact sequence in cohomology corresponding
to
$1{\rightarrow}{\mathbb Z}^n{\rightarrow}{\RB}^n{\rightarrow}{\mathbb 
T}^n{\rightarrow}1$, 
this implies that $H^2(\Gamma, {\mathbb Z}^n)=H^1(\Gamma, {\mathbb 
T}^n)$ which is finite 
and vanishes if we pass to a subgroup of finite index.  This allows 
us to lift the action of 
this subgroup of finite index, which is sufficient for our purposes.  
Note that the statement in the abstract follows from
the corollary, since Margulis' superrigidity theorem and the action 
on $H^1(M)$ being nontrivial imply that the 
action on ${\mathbb Z}^n$ is indeed given by the standard 
representation.

\begin{proof}
Since ${\rho}:\pi_1(M){\rightarrow}{\Lambda}$ is a surjection, we 
have a continuous map 
$f:M{\rightarrow}N/{\Lambda}$.
This follows since $N$ is contractible and $N/{\Lambda}$ is an 
Eilenberg-MacLane space for ${\Lambda}$ and hence
$M$ has a continuous map to $N/{\Lambda}$ inducing ${\rho}$ on 
fundamental groups which is canonical up to
homotopy.

We can lift $f$ to a $\Lambda$ equivariant map ${\tilde f}:{\tilde 
M}{\rightarrow}N$ where ${\tilde M}$ is
the cover of $M$ corresponding to $\ker{\rho}$.  

We consider the map $\tilde {\alpha}:{\Gamma}{\times}{\tilde 
M}{\rightarrow}N$ defined by 
$${\tilde {\alpha}}(\gamma, m)={\tilde 
f}({\gamma}m)({\gamma}f(m)){\inv}.$$
This is clearly a measure of the extent to which $\tilde f$ fails to 
be equivariant, and is defined since we have assumed the $\Gamma$ 
action lifts to $\tilde M$.  First we show that 
$\tilde {\alpha}$ descends to a map 
${\alpha}:\Gamma{\times}M{\rightarrow}N$.  To see this, let 
$m{\in}{\tilde M}$ and $m{\lambda}$ be any translate of $m$ where 
$\lambda{\in}{\Lambda}$ is viewed as a
deck transformation of $\tilde M$ over $M$.  It suffices to show that 
$\alpha(\gamma, m)={\alpha}(\gamma, m{\lambda})$, but
$${\alpha}(\gamma, m{\lambda})={\tilde 
f}({\gamma}m{\lambda})({\gamma}f(m{\lambda})){\inv}$$
$$={\tilde 
f}(({\gamma}m)({\gamma}{\lambda}))({\gamma}(f(m){\lambda})){\inv}$$
$$={\tilde 
f}({\gamma}m)({\gamma}{\lambda})(({\gamma}f(m)){\gamma}{\lambda})){\inv}$$

$$={\tilde 
f}({\gamma}m)({\gamma}{\lambda})({\gamma}{\lambda}){\inv}({\gamma}f(m)){\inv}$$

$$={\tilde f}({\gamma}m)({\gamma}f(m)){\inv}$$
since ${\tilde f}$ is ${\Lambda}$ equivariant and the action of 
$\Gamma$ on $\Lambda$ induced by the 
action on $\pi_1(M)$ is the same as the action of $\Gamma$ on 
$\Lambda<N$.

We now look at the map 
$\beta:{\Gamma}{\times}M{\rightarrow}{\Gamma}{\ltimes}N$ defined by 
${\beta}(\gamma, m)=(\gamma, \alpha(\gamma, m))$.  A simple 
computation verifies that $\beta$ is a 
cocycle over the $\Gamma$ action on $M$.  We can view $\beta$ as a 
cocycle into $G{\ltimes}N$ by the 
natural inclusion.  For now, we assume the action is ergodic. By 
results of Lewis and Zimmer, 
the algebraic hull, $L$, of this cocycle will be reductive 
with compact center.  Let $L_0<L$ be the connected component of the 
identity in $L$.  By passing to 
a finite ergodic extension of the action on $X=M{\times}{L/L_0}$ we 
have a  cocycle 
$\beta:{\Gamma}{\times}X{\rightarrow}G{\ltimes}N$ (still
called $\beta$) with algebraic hull $L_0$.  Note that $\beta(m,l)$ 
depends only on $m$.
Since any connected reductive subgroup of $G{\ltimes}N$ is conjugate 
to a subgroup of $G$, 
we can assume that $L_0<G$. 
This means that the cocycle on all of $X$ is cohomologous to one
taking values in $G$, in other words  $\beta(\gamma, 
x)=\phi(\gamma{x}){\inv}\delta(\gamma, x)\phi(x)$ 
where $\phi:X{\rightarrow}G{\ltimes}N$  is a measurable map and 
$\delta:{\Gamma}{\times}X{\rightarrow}G{\ltimes}N$ is a cocycle 
taking values 
entirely in $G$.  Write $\phi(x)=(\phi_1(x), \phi_2(x))=(\phi_1(x), 
1_N)(1_G, \phi_2(x))$.  
Since $\delta=(\delta_1, 1_N)$, computing the cocycle equivalence 
above in components yields that 
$\beta(\gamma, x)=(1_G, \phi_2(\gamma{x}){\inv}(\gamma, 1_N)(1_G, 
\phi_2(x))$.

We now show that ${\phi_2}$ is continuous.  
The argument follows exactly as in Lemma 6.5 of \cite{MQ}.  For the 
reader's convenience we repeat here the
case where $N={\mathbb R}^n$.  Essentially the idea is to use the 
fact that ${\mathbb R}^n$ is spanned by 
contracting directions for elements of $\Gamma$ and to show that 
along any contracting direction, $\phi_2$
can as the limit of iterated contractions of $\alpha$.

For $\gamma{\in}\Gamma$ let $E(\gamma)$ and $F(\gamma)$ be subspaces 
of ${\RB}^n$ that are the generalized
eigenspaces of $\gamma$ with eigenvalues of absolute value $>1$ and 
${\leq}1$ respectively.  Clearly 
${\RB}^n=E({\gamma}){\oplus}F(\gamma)$, and the assumption of weak 
hyperbolicity implies that ${\RB}^n$ is 
spanned by $\{E(\gamma)|\gamma{\in}{\Gamma}\}$.  We show continuity 
of $\phi_2$ by showing continuity of 
$\phi_2$
projected onto any $E(\gamma)$.  For any function 
$h:M{\rightarrow}{\RB}^n$ we write $h_{E(\gamma)}$ for the
composition of the function with projection on $E(\gamma)$.

Looking at what the cocycle condition on $\beta$ implies for $\alpha$ 
we see that 
$\phi_2(x)={\gamma}{\inv}{\phi_2}(\gamma{x})+{\gamma}{\alpha}(\gamma, 
m)$ where $x=(m,f){\in}X$.
Iterating this equality and projecting to $E({\gamma})$ gives
$${\phi_2}_{E(\gamma)}(x)={\sum_{i=1}^n(\gamma^i){\inv}|_{\eg}{\alpha_{\eg}}(\gamma, 
{\gamma}^{i-1}m)}
+ {\gamma}^n{\inv}|_{\eg}{\phi_2}_{\eg}({\gamma}^nx).$$
Since the eigenvalues of ${\gamma}{\inv}|_{\eg}$ all have absolute 
value $<1$,
on a set of full measure  
${\gamma}^n{\inv}|_{\eg}{\phi_2}_{\eg}({\gamma}^nx){\rightarrow}0$
 as $n{\rightarrow}{\infty}$.  So we have:
\begin{equation}\tag{$*$}
{\phi_2}_{E(\gamma)}(x)={\sum_{i=1}^{\infty}(\gamma^i){\inv}|_{\eg}{\alpha_{\eg}}(\gamma, 
{\gamma}^{i-1}m)}
\end{equation}
which converges uniformly since $\alpha(\gamma, -)$ is continuous and 
bounded function on $M$.  This shows 
both that $\phi_2$ is continuous and is a function on $M$ that is 
independent of the finite ergodic extension
$X$.

If the action is not ergodic, we simply carry out the analysis above 
on each ergodic component.  We will
get a function $\phi_2$ as above for each component.  Since $(*)$ 
above shows how to compute $\phi_2$
explicitly on any ergodic component, we see that $\phi_2$ is a well 
defined continuous function from $X$ to
${\RB}^n$.  From $(*)$ it is clear that $\phi_{2}$ descends to a 
continuous function $M$ to ${\RB}^{n}$. When ${\RB}^n$ is replace by a 
more general simply connected nilpotent group, the analysis becomes
more complicated but follows exactly as in \cite{MQ}.

Now $\beta(\gamma, m)=(\gamma, {\tilde f}(\gamma{m})({\gamma}{\tilde 
f}(m))\inv)$ where we are actually
choosing some lift of the point $m$ to $\tilde M$.  Substituting this 
in above and computing the $N$ 
factor gives ${\tilde f}({\gamma}m)(\gamma{\tilde 
f}(m))\inv={\phi_2}({\gamma}m){\inv}(\phi_2(m))$. 
 We can lift
$\phi_2$ to a map from $\tilde M$ to $N$, and then rearranging the 
last expression gives 
${\tilde {\phi_2}}({\gamma}m){\tilde f}({\gamma}m)={\gamma}{\tilde 
{\phi_2}}(m){\gamma}{\tilde f}(m)=\gamma({\tilde {\phi_2}}{\tilde 
f})$.  
This shows that the map $(\tilde {\phi_2})(\tilde f):{\tilde 
M}{\rightarrow}N$ is $\Gamma$ 
equivariant.  Since ${\tilde {\phi_2}}$ is $\Lambda$ invariant, 
$\tilde f$ is $\Lambda$ equivariant and 
$\Lambda$ acts on the right on $\tilde M$, we see that the map 
${\phi_2}f:M{\rightarrow}N/{\Lambda}$ is also 
$\Gamma$ equivariant.  Note that $\phi_2:M{\rightarrow}N$ is a map 
into a contractible space and so isotopically
 trivial.  This implies that $\phi_2f$ is in the same isotopy class 
as $f$ and therefore that 
${\phi_2}f$ induces the map $\rho$ on fundamental group.

\end{proof}

\section{{\bf Examples}}
\label{section:a&e}

For the sake of clarity we discuss the case of Corollary 
\ref{corollary:torus} only, although much can easily be 
generalized to the general case of Theorem \ref{theorem:main}.  
Throughout $\Gamma$ will refer to a finite index subgroup of 
$SL_{n}{\mathbb Z}$.

We start with some non-trivial examples to which our theorem applies.
These examples of exotic $SL_{n}({\mathbb Z})$ actions are due to 
Katok and Lewis (\cite{KL}).  They are constructed from the 
standard action on ${\mathbb T}^{n}$ by blowing up the fixed point, 
so that it becomes a copy of ${\mathbb RP}^{n-1}$ with the standard 
action of $SL_{n}{\mathbb Z}$.  What is not at all obvious that the 
resulting manifold has an invariant analytic structure and volume form.  
Indeed, in order to make the action preserve a volume form one must 
use a differential structure which is not the obvious one.  The 
continuous map to the torus,  collapsing the ${\mathbb RP}^{n-1}$, is 
not smooth with respect to this new smooth structure.  It is immediate 
from the proof of Theorem \ref{theorem:main} that the map is 
unique, since it is given explicitly by (*) in terms of the action.  This 
shows that the regularity of the semi-conjugacy in Theorem \ref{theorem:main} 
cannot be improved, even when the action is analytic.

There are further examples of exotic actions of $SL_n({\mathbb Z})$, 
which were originally constructed by Weinberger.  Take ${\mathbb 
T}^{n}$ and remove some finite invariant set.  The resulting manifold 
can be compactified to a manifold with 
boundary by adding the spheres in the tangent bundle at each point.  
The action of $\Gamma$ on the boundary is the standard linear
actions.  These actions extend over $n$-balls - just think of the 
ball 
as $\RB^{n}$ with the linear $\Gamma$ action as the action on the 
sphere at 
infinity.  Gluing in these balls gives a closed manifold with a 
$\Gamma$ action.  The underlying manifold is still the torus, but the 
action is different.
There is no invariant measure on the whole torus, but there is on the 
complement of the disks.  Thus our theorem says there is a continuous 
map
from this complement to the torus.  Indeed the map which simply 
collapses 
the balls we glued in back to points is continuous and equivariant. 

It is possible to combine Weinberger's construction with some surgery 
to produce a new class of examples.  Take $(X,\partial X)$ any 
compact manifold with boundary, let $N$ be Weinberger's example 
cross a $\partial X$.  Inside of $N$ is a ball cross $\partial X$.  
Remove the interior of this, and glue in $X$ cross an $n-1$ sphere in 
its place to get $M$.  If $\pi_{1}(X)$ is trivial, $\pi_1(M)={\mathbb 
Z}^{n}$.  As a specific example, taking $X$ to be the $m$ ball gives 
an example of a $\Gamma$-manifold with $\pi_{1}={\mathbb Z}^{n}$ but 
which is definitely not, even non-equivariantly, a bundle over 
${\mathbb T}^{n}$.

In the above examples the map to the torus is always well behaved
off of a submanifold which is collapsed.  The following 
lemma, in some ways a converse to the main theorem,  shows that this
collapsing must occur in general.

\begin{lemma} Let $Z$ be any compact, connected, space with a 
$\Gamma$ action.  
Any non-constant equivariant map to the torus surjects $\pi_{1}(Z)$ 
onto a finite index subgroup of ${\mathbb Z}^{n}$.
\begin{proof}

Suppose not.  Then, since the image of $\pi_{1}(Z)$ is a $\Gamma$ 
invariant infinite index subgroup of ${\mathbb Z}^{n}$, it must be 
trivial.  Thus we aim to show that any null homotopic map is constant.
Let $f:Z \to {\mathbb T}^{n}$ be a null homotopic equivariant map.

The image of $f$ in ${\mathbb T}^{n}$ is closed, invariant, and 
connected.  Hence, since we assume the map non-constant, it must be
surjective.  In particular, some $z_{0}$ in $M$ maps to $0$ in 
${\mathbb 
T}^{n}$.  Lift $f$ to a map $F:Z \to \RB^{n}$ such that $F(z_{0})=0$.

Fix a $\gamma \in \Gamma$.  Since $F$ covers an equivariant map, we 
know that $$\tau_{\gamma}(z)=F(\gamma z)-\gamma F(z)$$ takes values 
in ${\mathbb Z}^{n}$.  Since it is continuous in $z$, this implies 
it is constant in $z$.  We define $\tau_{\gamma}$ to be the common 
value.
It is easy to see that $\tau : \Gamma \to {\mathbb Z}^{n}$ is a 
1-cocycle, 
in other words $$\tau_{\gamma \sigma}= \tau_{\gamma} + \gamma 
\tau_{\sigma}.$$

We can evaluate $\tau_{\gamma}$ at $z_{0}$, which yields 
$\tau_{\gamma}=F(\gamma z_{0})$.  Since the image of $F$ is compact, 
this shows that $\tau_{\gamma}$ is bounded indepependant of 
$\gamma$.  The cocycle identity then shows that $\gamma 
\tau_{\sigma}$ 
is bounded independantly of $\gamma$ and $\sigma$.  In particular, 
$\tau_{\sigma}$ has bounded $\Gamma$ orbit and so we have 
$\tau_{\sigma}=0$ for all $\sigma$.  This means precisely that $F$ is 
equivariant.  This finishes the proof as $F(Z)$ is then a bounded, 
invariant set in $\RB^{n}$, and thus is $\{0\}$.
 
\end{proof}
\end{lemma}

In all of these examples, the map to the torus is nice everywhere
except the pre-image of a finite invariant set.  We believe that will 
always be the case.  Note that if the map had any regularity, then by 
Sard's theorem the critical values would be measure zero.  Since they are 
closed and invariant this would limit them to a finite invariant set, 
and thus we would know that off this finite set our manifold is a 
bundle over the torus.  We know from the Katok-Lewis examples that 
the map need not be $C^{1}$, even for analytic actions.  The map in 
that case is, however, still analytic off a lower dimensional 
submanifold.

\begin{question}  If, in the statement of the main theorem, one 
assumes the action to have some regularity, is the map also regular 
away from the pre-image of a finite invariant set in the torus?
\end{question}

Even if regularity does not hold, one can still hope the map must be 
``taut'' in some sense.  The last lemma is one example of the sort of
substitute for regularity.  At the very least one would like to be able
to rule out space filling curves in this context, so as to prove:

\begin{conjecture} Any compact manifold with $\pi_{1}={\mathbb 
Z}^{n}$ 
and a $\Gamma$ action which induces the standard action on $\pi_{1}$ 
must be of dimension at least $n$.
\end{conjecture}

  Here we have no assumption of an invariant measure.  In all the 
examples there is an invariant measure, at least on a large open 
set.  Is this always the case?

\begin{question} If $SL_{n}({\mathbb Z})$ acts on a compact manifold
with fundamental group ${\mathbb Z}^{n}$, inducing the standard action
on $\pi_{1}$, is there always an invariant measure? Is there always 
an invariant measure whose support contains an open set?
\end{question}

  One is tempted to view the torus with the standard action as some kind 
of equivariant classifying space for $\Gamma$ actions with $\pi_{1}={\mathbb 
Z}^{n}$, and to view our theorem as proving the existence of a 
classifying map  One cannot hope for general results of this type - 
the use of superrigidity is not just an artifact of the method. Even
when there is a clear candidate classifying space, and the action
there is hyperbolic, the analog of our theorem need not hold.

  Consider the action of $SL_{2}({\mathbb Z})$ on ${\mathbb 
T}^{2}$.   Since $SL_{2}$ acts hyperbolically on the torus, one might 
expect our theorem to cover this case.  It does not: take $\Gamma$ any 
torsion free subgroup of finite index in $SL_{2}({\mathbb Z})$.  Such 
a $\Gamma$ is free.  Starting with the standard action 
of $\Gamma$ on ${\mathbb T}^{2}$, conjugate the action of 
one of the free generators by a homeomorphism homotopic to the 
identity, and leave the action of the remaining generators unchanged.  
Since $\Gamma$ is free, this still generates an action of $\Gamma$. 
There is no equivariant map of this torus to the standard one.  This 
follows from the uniqueness of the map conjugating a single Anosov 
homeomorphism to a linear Anosov automorphism, which follows from
the same reasoning that shows the map in Theorem \ref{corollary:torus}
is unique (or see \cite{KH}).

\medskip
\noindent
David Fisher\\
Department of Mathematics\\
Yale University\\
P.O. Box 208283\\
New Haven, CT 06520-8283\\
E-mail: david.fisher@yale.edu

\medskip
\noindent
Kevin Whyte\\
Department of Mathematics\\
University of Chicago\\
5734 S. University Ave\\
Chicago, Il 60637\\
E-mail: kwhyte@math.uchicago.edu\\

\end{document}